\documentclass[letterpaper, 10 pt, conference]{ieeeconf}    \IEEEoverridecommandlockouts    

\usepackage{graphicx}
\usepackage{epsfig} 
\usepackage{times} 
\usepackage{amsmath}
\usepackage{amssymb}  
\usepackage{color}
\usepackage{amsfonts}
\usepackage{subfigure}
\usepackage{multirow}
\usepackage{multicol}
\usepackage{siunitx}
\usepackage{overpic}
\usepackage{mathrsfs}  

\usepackage{booktabs}
\usepackage{threeparttable}

\usepackage{epstopdf}
\usepackage{xspace}

\begin{document}


\title{\LARGE \bf Highway Traffic State Estimation with Mixed \\ Connected and Conventional Vehicles}

\author{Nikolaos Bekiaris-Liberis, Claudio Roncoli, and Markos Papageorgiou
\thanks{N. Bekiaris-Liberis, Claudio Roncoli, and Markos Papageorgiou are with the Department of Production Engineering and Management, Technical University of Crete, Chania, 73100, Greece. Email addresses: \texttt{nikos.bekiaris@gmail.com, croncoli@dssl.tuc.gr, {\rm and} markos@dssl.tuc.gr}.} }

\maketitle
\thispagestyle{empty}

\begin{abstract}
A macroscopic model-based approach for estimation of the traffic state, specifically of the (total) density and flow of vehicles, is developed for the case of ``mixed" traffic, i.e., traffic comprising both ordinary and connected vehicles. The development relies on the following realistic assumptions: (i) The density and flow of connected vehicles are known at the (local or central) traffic monitoring and control unit on the basis of their regularly reported positions; and (ii) the average speed of conventional vehicles is roughly equal to the average speed of connected vehicles. Thus, complete traffic state estimation (for arbitrarily selected segments in the network) may be achieved by merely estimating the percentage of connected vehicles with respect to the total number of vehicles. A model is derived, which describes the dynamics of the percentage of connected vehicles, utilizing only well-known conservation law equations that describe the dynamics of the density of connected vehicles and of the total density of all vehicles. Based on this model, which is a linear time-varying system, an estimation algorithm for the percentage of connected vehicles is developed employing a Kalman filter. The estimation methodology is validated through simulations using a second-order macroscopic traffic flow model as ground truth for the traffic state. The approach calls for a minimum of spot sensor-based total flow measurements {according to a variety of possible location configurations.}
\end{abstract}
\section{Introduction}
A number of novel Vehicle Automation and Communication Systems (VACS) have already been introduced, and many more are expected to be introduced in the next years. These systems are mainly aimed to improve driving safety and convenience, but are also believed to have great potential in mitigating traffic congestion, if appropriately exploited for innovative traffic management and control \cite{diakaki}. To attain related traffic flow efficiency improvements on highways, it is of paramount importance to develop novel methodologies for modeling, estimation and control of traffic in presence of VACS. Several papers are providing useful results related to modeling and control of traffic flow in presence of VACS, employing either microscopic or macroscopic approaches, see, for example, {\cite{bose1}}, \cite{bose2}, \cite{davis1}, \cite{ge1}, \cite{kesting}, \cite{ngoduy4}, \cite{rajamani}, \cite{rao}, \cite{roncoli0}, \cite{roncoli00}, \cite{roncoli1}, \cite{shladover}, \cite{vanarem1}, \cite{varaya1}, \cite{wang van arem1}, \cite{yi}.

The availability of reliable real-time measurements or estimates of the traffic state is a prerequisite for successful highway traffic control. In conventional traffic, the necessary measurements are provided by spot sensors (based on a variety of possible technologies), which are placed at specific highway locations.  If the sensor density is sufficiently high (e.g., every 500 m), then the collected measurements are usually sufficient for traffic surveillance and control; else, appropriate estimation schemes need to be employed in order to produce traffic state estimates at the required space resolution (typically 500 m); see, for instance, \cite{alvarez}, \cite{heygi1}, \cite{miha}, \cite{munoz}, \cite{papa new}, among many other works addressing highway traffic estimation by use of conventional detector data. However, the implementation and maintenance of road-side detectors entail considerable cost; hence various research works attempt to exploit different, less costly data sources, such as mobile phone, or GPS (Global Positioning System), or even vehicle speed data for travel time or highway state estimation; see, e.g., \cite{atsarita}, \cite{fabri}, \cite{deng}, \cite{bayen1}, \cite{bayen2}, \cite{ou}, \cite{rahmani1}, \cite{seo}, \cite{treiber}, \cite{work1}, \cite{yuan2}; employing various kinds of traffic or statistic models.

In fact, with the introduction of VACS of various kinds, an increasing number of vehicles become ``connected", i.e., enabled to send (and receive) real-time information to a local or central monitoring and control unit (MCU). Thus, connected vehicles may communicate their position, speed and other relevant information, i.e., they can act as mobile sensors. This will potentially allow for a sensible reduction (and, potentially, elimination) of the necessary number of spot sensors, which would lead to sensible reduction of the purchase, installation, and maintenance cost for traffic monitoring. This paper concerns the development of reliable and robust traffic state estimation methods and tools, which exploit information provided by connected vehicles and reduces the need for spot sensor measurements under all penetration rates of connected vehicles, i.e., for a mixed traffic flow that includes both conventional and connected vehicles.

Specifically, we address the problem of estimating the (total) density and flow of vehicles in highway segments of arbitrary length (typically around 500 m) in presence of connected vehicles. The developments rely on the following realistic assumptions:
\begin{itemize}
\item The density and flow of connected vehicles may be readily obtained at the local or central MCU on the basis of their regularly reported positions. 
\item The average speed of conventional vehicles is roughly equal to the average speed of connected vehicles. This assumption relies on the fact that, even at very low densities, there is no reason for connected vehicles to feature a systematically different mean speed than conventional vehicles; while at higher densities, the assumption is further reinforced due to increasing difficulty of overtaking.
\end{itemize}
As a consequence of these assumptions, complete traffic state estimation (of the total density and flow in arbitrarily selected segments in the highway) may be achieved by merely estimating the percentage of connected vehicles with respect to the total number of vehicles. For the latter, a minimum amount of conventional measurements of traffic volumes, e.g., at all highway entries and exits, is also required. Thus, the problem of traffic estimation is recast in the problem of estimating the percentage of connected vehicles at the selected highway segments.

    In more technical terms, we derive a linear time-varying {model}, which describes the dynamics of the percentage, utilizing merely the {(time-discrete)} conservation law equations for the density of connected vehicles and for the total density of vehicles (no traffic modelling of speed, such as the fundamental diagram, is required). We show that the system is observable and employ a Kalman filter for the estimation of the percentage of connected vehicles. We demonstrate our estimation design with a numerical example employing a second-order macroscopic traffic flow model as ground truth for the traffic state dynamics. 

     Section \ref{sec model} derives a linear time-varying system that describes the dynamics of the percentage of connected vehicles at each segment of a highway. Section \ref{seccal1} studies the observability properties of the system; while Section \ref{seccal2} employs a Kalman filter for the estimation of the percentage of connected vehicles on the highway. Section \ref{seccal3} demonstrates the estimation design with a second-order macroscopic model as ground truth. Section \ref{uns} extends the approach to the case of unmeasured total flow at off-ramps. Section \ref{sec conclusions} summarizes the conclusions and outlines related ongoing and future work.

\section{Model Derivation for the Percentage of connected Vehicles}
\label{sec model}
We consider the following discrete-time equations that describe the dynamics of the total density $\rho$ of the vehicles on a highway and the density $\rho^{\rm a}$ of the connected vehicles {(see, for example, \cite{papageorge11}; see also the upper part of Fig. \ref{fig1})}
\begin{eqnarray}
\rho_i(k+1)&=&\rho_i(k)+\frac{T}{\Delta_i}\left(q_{i-1}(k)-q_i(k)+r_i(k)\right.\nonumber\\
&&\left.-s_i(k)\right)\label{eqrho}\\
\rho^{\rm a}_i(k+1)&=&\rho^{\rm a}_i(k)+\frac{T}{\Delta_i}\left(q^{\rm a}_{i-1}(k)-q^{\rm a}_i(k)+r^{\rm a}_i(k)\right.\nonumber\\
&&\left.-s^{\rm a}_i(k)\right),\label{eqrhoa}
\end{eqnarray}
where $i=1,\ldots,N$ is the index of the specific segment at the highway, $N$ being the number of discrete cells on the {highway;} for all traffic variables, we denote by index sub-$i$ its value at the segment $i$ of the highway; $q_i$ and $q^{\rm a}_i$ are the total flow and the flow of the connected vehicles, respectively, at segment $i$; $T$ is the time-discretization step, $\Delta_i$ is the length of the discrete segments of the highway, and $k=0,1,\ldots$ is the discrete time index. The variables $r_i$ and $s_i$ denote the inflow and outflow of vehicles at on-ramps and off-ramps, respectively, at segment $i$, whereas $r_i^{\rm a}$ and $s_i^{\rm a}$ are the corresponding inflow and outflow of connected vehicles. Define the inverse of the percentage of the connected vehicles at segment $i$ of the highway as $\bar{p}_i$, i.e.,
\begin{eqnarray}
\bar{p}_i&=&\frac{\rho_i}{\rho^{\rm a}_i}\label{eq1}.
\end{eqnarray}
{Assuming} that the average speed of conventional vehicles at a segment $i$ equals the average speed of connected vehicles in the same segment, namely $v_i$, one can conclude that the following holds
\begin{eqnarray}
\bar{p}_i=\frac{\rho_i}{\rho^{\rm a}_i}=\frac{q_i}{q^{\rm a}_i},\label{aa}
\end{eqnarray} 
where we used the known relations
\begin{eqnarray}
q_i&=&\rho_iv_i\label{flowtotal}\\
q^{\rm a}_i&=&\rho^{\rm a}_iv_i.\label{flowvacs}
\end{eqnarray}
Using (\ref{eqrho}), (\ref{eqrhoa}), and (\ref{aa}) we get from (\ref{eq1}) that
\setlength{\arraycolsep}{0pt} \begin{eqnarray}
 \bar{p}_i(k+1)&=&\frac{\left(\rho^{\rm a}_i(k)-\frac{T}{\Delta_i}q^{\rm a}_i(k)\right)\bar{p}_i(k)+\frac{T}{\Delta_i}q^{\rm a}_{i-1}(k)\bar{p}_{i-1}(k)}{g_i^{\rm a}(k)}\nonumber\\
 &&+\frac{T}{\Delta_i}\frac{\left(r_i(k)-s_i(k)\right)}{g_i^{\rm a}(k)}\label{per1}\\
 g_i^{\rm a}(k)&=&\rho^{\rm a}_i(k)+\frac{T}{\Delta_i}\left( q^{\rm a}_{i-1}(k)-q^{\rm a}_i(k)+r^{\rm a}_i(k)\right.\nonumber\\
 &&\left.-s^{\rm a}_i(k)\right),\label{defcru}
 \end{eqnarray}\setlength{\arraycolsep}{5pt}$i=1,\ldots,N$. Defining the state 
\begin{eqnarray}
{x}=\left(\bar{p}_1,\ldots,\bar{p}_N\right)^T,
\end{eqnarray}
 we re-write (\ref{per1}) as 
\begin{eqnarray}
{x}(k+1)&=&{A}(k){x}(k)+{B}(k){u}(k)\label{12}\\
{y}(k)&=&{C}{x}(k),\label{120}
\end{eqnarray}
where 
\setlength{\arraycolsep}{4pt}\begin{eqnarray}
{A}(k)&=&\left\{\begin{array}{lll}{a}_{ij}=\frac{T}{\Delta_i}\frac{q_{i-1}^{\rm a}(k)}{g_i^{\rm a}(k)},&\mbox{if $i-j=1$}\\&\mbox{and $i\geq2$}\\{a}_{ij}=\frac{\rho^{\rm a}_i(k)-\frac{T}{\Delta_i}q^{\rm a}_i(k)}{g_i^{\rm a}(k)},&\mbox{if $i=j$}\\{a}_{ij}=0,&\mbox{otherwise}\end{array}\right\}\label{adef}\\
{B}(k)&=&\left\{\begin{array}{lll}{b}_{ij}=\frac{T}{\Delta_i}\frac{1}{g_1^{\rm a}(k)},&\mbox{if $i=1$}\\&\mbox{and $j=1,2$}\\{b}_{ij}=\frac{T}{\Delta_i}\frac{1}{g_i^{\rm a}(k)},&\mbox{if $j-i=1$}\\{b}_{ij}=0,&\mbox{otherwise}\end{array}\right\}\\
{u}(k)&=&\left[\begin{array}{c}q_0(k)\\r_1(k)-s_1(k)\\\vdots\\r_N(k)-s_N(k)\end{array}\right]\label{dubu}\\
{C}&=&{\left[\begin{array}{cccc}0&\ldots&0&1\end{array}\right]},\label{16}
\end{eqnarray}\setlength{\arraycolsep}{5pt}$g_i^{\rm a}$, $i=1,\ldots,N$, is defined in (\ref{defcru}), $A\in\mathbb{R}^{N\times N}$, $B\in\mathbb{R}^{N\times (N+1)}$, and $q_0$ denotes the total flow of vehicles at the entry of the highway 
and acts as an input to system (\ref{12}), along with the variables $r_i$ and $s_i$; while $r^{\rm a}_i$, $s^{\rm a}_i$, $\rho_i^{\rm a}$, and $q_i^{\rm a}$ are viewed as time-varying parameters of system (\ref{12}). Finally, the variable $\bar{p}_N$ at the exit of the highway is viewed as the output of the system and may be obtained via 
\begin{eqnarray}
\bar{p}_N=\frac{q_N}{q^{\rm a}_N},
\end{eqnarray}
using total flow measurements $q_N$ at the highway exit. 

Before studying the observability of system (\ref{12})--(\ref{16}), we summarize the assumptions that guarantee that the matrix ${A}$ is known, and that the input ${u}$ and output ${y}$ are measured.
\begin{itemize}
\item The average speed of the connected vehicles at a segment of the highway equals the average speed of all vehicles at the same segment, i.e., $v_i^{\rm a}=v_i$.
\item {The segment flows and densities} of connected vehicles, $q_i^{\rm a}$, $i=0,\ldots,N$, and $\rho_i^{\rm a}$, $i=1,\ldots,N$, respectively, as well as the flows of connected vehicles at on-ramps and off-ramps, $r_i^{\rm a}$ and $s_i^{\rm a}$, $i=1,\ldots,N$, respectively, may be obtained from regularly received messages by the connected vehicles.
\item The total flow of vehicles at the entry and exit of the highway, $q_0$ and $q_N$, respectively, are measured via conventional detectors.
\item The total flow of vehicles at on-ramps and off-ramps, $r_i$ and $s_i$, $i=1,\ldots, N$, respectively, are measured via conventional detectors. 
\end{itemize}
The above formulation may be modified in a couple of respects:
\begin{itemize}
\item Different total flow measurement configurations may be employed; for example, additional mainstream total flow measurements (using conventional detectors) may be considered to replace a corresponding number of total flows at on-ramps or off-ramps, without affecting the observability of the system.
\item In case more mainstream total flow measurements are actually employed (in place of total flow measurements at on- and off-ramps), they may also be considered as output variables in (\ref{120}) to potentially increase the filterÕs responsiveness.
\end{itemize}
These issues are currently in the course of investigation.


Note that later on, in Section \ref{uns}, we remove, under certain conditions, the assumption that the total off-ramp flows $s_i$, $i=1,\ldots,N$ are measurable. 

\section{Percentage Estimation Using a Kalman Filter}
\label{seccal}
\subsection{Observability of the System}
\label{seccal1}
System (\ref{12}) is viewed as a linear time-varying system. As it is stated in Section \ref{sec model}, it is assumed that the quantities ${q}_0$, $\bar{p}_N$, $q^{\rm a}_i$, $\rho^{\rm a}_i$, $r^{\rm a}_i$, $s^{\rm a}_i$, $r_i$, and $s_i$, for all $i$, are available, which implies that the matrices ${A}$ and ${B}$, as well as the input ${u}$ in (\ref{12}) may be calculated in real time. We show next that system (\ref{12})--(\ref{16}) is observable at $k=k_0+N-1$, for any initial time $k_0\geq0$. We construct the observability matrix
\begin{eqnarray}
{{O}(k_0,k_0+{N})}=\left[\begin{array}{cc}{C}\\{C}{A}(k_0)\\{C}{A}(k_0+1){A}(k_0)\\\vdots\\{C}{A}(k_0+{N}-2)\cdots {A}(k_0)\end{array}\right].
\end{eqnarray}
{Since ${O}$} is square, the system is observable at $k=k_0+N-1$ if $\det({O})\neq0$. Since from (\ref{adef}) it is evident that ${A}$ is a lower triangular matrix, it follows from (\ref{16}) that ${O}$ is an anti-lower triangular matrix, namely, a matrix with zero elements above the anti-diagonal. Therefore, relation $\det({O})\neq0$ holds if the anti-diagonal elements of ${O}$ are non-zero. The anti-diagonal elements of ${O}$ are given by $1, {a}_{NN-1}(k_0), {a}_{NN-1}(k_0+1){a}_{N-1N-2}(k_0),\ldots,{a}_{NN-1}(k_0+N-2)\cdots {a}_{21}(k_0)$. Since $q_i^{\rm a}$, $\rho_i^{\rm a}$, $i=1,\ldots,N$, are lower and upper bounded (and positive), it follows from (\ref{adef}) that ${a}_{ij}(k)$, for all $k=k_0,\ldots,k_0+N-2$ and any $k_0\geq0$, and for all $i,j$ such that $i-j=1$ and $i\geq2$, are lower and upper bounded (and positive). Therefore, the matrix ${O}$ is invertible, and hence, system (\ref{12})--(\ref{16}) is observable at $k=k_0+N-1$. Note that the measurement of $\bar{p}_N$ (or, equivalently, of $q_N$), rather than any other intermediate percentage, is necessary for system (\ref{12})--(\ref{16}) to be observable. To see this note that if ${C}=\left\{\begin{array}{cc}{c}_{ij}=1,&\mbox{if $i=1$ and $j=J$}\\{c}_{ij}=0,&\mbox{otherwise}\end{array}\right\}$ with $J<N$, then the $J+1,\ldots,N$ columns of {${O}(k_0,k_0+\bar{N})$ are zero for all $k_0\!\geq\!0$ and $\bar{N}\!\geq\! N$}. Thus, the system cannot be observable. In other words, a fixed flow sensor should necessarily be placed at the last segment of the highway in order to guarantee percentage observability based on model (\ref{12})--(\ref{16}).

\subsection{Kalman Filter}
\label{seccal2}
\begin{figure}
\centering
\includegraphics[width=\linewidth]{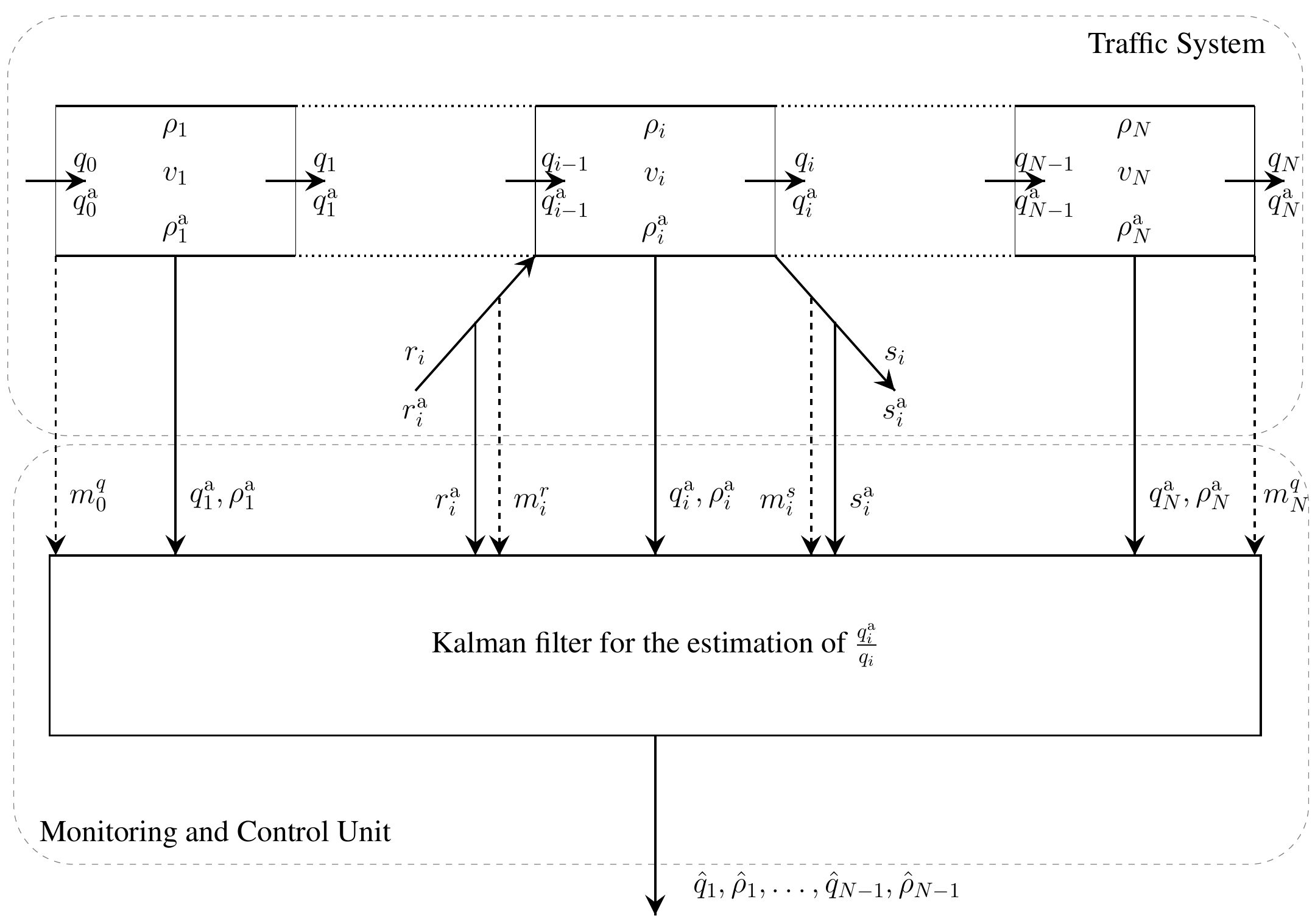}
\caption{The traffic system under consideration and the Kalman filter implemented at the MCU. The data used to operate the Kalman filter are either coming from connected vehicles (solid lines) or fixed sensors (dashed lines). The variable $m_i^w$ denotes the measurement of quantity $w$ at segment $i$, which might be different than the actual quantity $w$, due to, for example, the presence of measurement noise.}
\label{fig1}
\end{figure}

We implement a Kalman filter for the estimation of the percentage of connected vehicles on a highway (see Fig. \ref{fig1}). Defining $\hat{{x}}=\left(\hat{\bar{p}}_1,\ldots,\hat{\bar{p}}_N\right)^T$, the equations for the Kalman filter are given by (see, for example, \cite{anger})
\begin{eqnarray}
\hat{{x}}(k+1)\!&=&\!{A}(k)\hat{{x}}(k)+{B}(k){u}(k)\nonumber\\
&&+{{A}}(k){K}(k)\left({{z}}(k)-{C}\hat{{x}}(k)\right)\label{123}\\
{K}(k)\!&=&\!{P}(k){C}^T\left({C}{P}(k){C}^T+{R}\right)^{-1}\\
{P}(k+1)\!&=&\!{A}(k)\left(I-{K}(k){C}\right){P}(k){A}(k)^T+{Q}\!,\label{kal1}
\end{eqnarray}
where ${z}$ is a noisy version of the measurement ${y}$, ${R}>0$ and ${Q}={Q}^T>0$ are tuning parameters. Note that, in the ideal case in which there is additive, zero-mean Gaussian white noise in the output and state equation (\ref{12}) and (\ref{120}), respectively, these matrices represent the (ideally known) covariance matrices of the measurement and process noise, respectively. Since the system equations here are relatively complex, some tuning of the matrices may be necessary for best estimation results. 
The initial conditions of the estimator (\ref{123})--(\ref{kal1}) are given by
\begin{eqnarray}
\hat{{x}}(k_0)&=&\mu\\
{P}(k_0)&=&H,\label{1234}
\end{eqnarray}
where $\mu$ and $H=H^T>0$ are the initial conditions of the estimator (\ref{123})--(\ref{kal1}), which, in the ideal case in which ${x}(k_0)$ is a Gaussian random variable, represent the mean and auto covariance matrix of ${x}(k_0)$, respectively.

The Kalman filter (\ref{123})--(\ref{1234}) delivers estimates of the inverse percentages $\hat{\bar{p}}_i$; using (\ref{aa}) and the available data for $q_i^{\rm a}$, $\rho_i^{\rm a}$, we can obtain estimates for all segment (total) flows and densities $\hat{q}_i$, $\hat{\rho}_i$ as indicated at the output of the Kalman filter in Fig. \ref{fig1}.



\subsection{Evaluation of the Performance of the Estimator Based on a METANET Model as Ground Truth}
\label{seccal3}
For preliminary assessment of the developed estimation scheme, we test in this section the performance of the Kalman filter employing the second-order METANET model \cite{papageorge11} (i.e., a model in which the average speed of the vehicles at the highway has its own dynamics) as ground truth. We employ equations (\ref{eqrho}) and (\ref{eqrhoa}) for the total density of the vehicles and the density of connected vehicles, respectively, together with relations (\ref{flowtotal}) and (\ref{flowvacs}) for the total flow and the flow of connected vehicles, respectively. The equation for the average speed at segment $i$ is given by
\setlength{\arraycolsep}{0pt}\begin{eqnarray}
v_i(k+1)&=&v_i(k)+\frac{T}{\tau}\left(V\left(\rho_i(k)\right)-v_i(k)\right)+\frac{T}{\Delta_i}v_i(k)\nonumber\\
&&\times\left(v_{i-1}(k)-v_i(k)\right)-\frac{\nu T}{\tau\Delta_i}\frac{\rho_{i+1}(k)-\rho_i(k)}{\rho_i(k)+\kappa}\nonumber\\
&&-\frac{\delta T}{\Delta_i}\frac{r_i(k)v_i(k)}{\rho_i(k)+\kappa},\quad i=1,\ldots,N,\label{average speed}
\end{eqnarray}\setlength{\arraycolsep}{5pt}with $v_0=v_1$ and $\rho_N=\rho_{N+1}$, where the nominal average speed $V$ is given by
\begin{eqnarray}
V\left(\rho\right)=v_{\rm f}e^{-\frac{1}{\alpha}\left(\frac{\rho}{\rho_{\rm cr}}\right)^{\alpha}},\label{fund}
\end{eqnarray}
and $\tau$, $\nu$, $\kappa$, $\delta$, $v_{\rm f}$, $\rho_{\rm cr}$, and $\alpha$ are positive model parameters. In particular, $v_{\rm f}$ denotes the free
speed, $\rho_{\rm cr}$ the critical density, and $\alpha$ the exponent of the stationary speed equation (\ref{fund}). The model parameters, which are taken from \cite{papageorge1}, are shown in Table \ref{table1}.
 \begin{table}[t]
\caption{Parameters of the model (\ref{eqrho}), (\ref{eqrhoa}), (\ref{flowtotal}), (\ref{flowvacs}), (\ref{average speed}), and (\ref{fund}).}
\begin{center}
\begin{tabular}{cl}
\hline\hline
Model parameter & Value\\
\hline
$T$&$\frac{10}{3600}$ (h)\\[2mm]
$\Delta_i$&$\frac{500}{1000}$ ({km})\\[2mm]
$\tau$&$\frac{20}{3600}$ (h) \\[2mm]
$\nu$&$35$ $\left(\frac{\textrm{km}^{\textrm{2}}}{\textrm{h}}\right)$ \\[2mm]
 $\kappa$&$13\left(\frac{\textrm{veh}}{\textrm{km}}\right)$\\[2mm]
 $\delta$&$1.4$\\[2mm]
 $v_{\rm f}$&$120\left(\frac{\textrm{km}}{\textrm{h}}\right)$\\[2mm]
 $\rho_{\rm cr}$&$33.5\left(\frac{\textrm{veh}}{\textrm{km}}\right)$\\[2mm]
 $\alpha$&$1.4324$\\[2mm]
  $N$&$20$\\[2mm]
\hline
\end{tabular}
\label{table1}
\end{center}
\end{table}%

From the model parameters (\ref{adef})--(\ref{16}) and the Kalman filter (\ref{123})--(\ref{kal1}) it is evident that the estimator utilizes measurements stemming from connected vehicles reports, namely, $q^{\rm a}_i$, $\rho^{\rm a}_i$, $r^{\rm a}_i$, $s^{\rm a}_i$, for all $i$. Although the calculation of these variables from connected vehicle data is likely to be associated with error or noise, we assume, for this preliminary assessment, that they are accurate measurements. In contrast, the measurements of the total flow of the vehicles at the entry and exit of the highway are subject to additive measurement noise, say, $\gamma_0^q\sim N(0,D_q^2)$ and $\gamma_N^q\sim N(0,D_q^2)$, respectively. Furthermore, the measurements of the total flow at the on-ramps and off-ramps might be subject to additive measurement noise say $\gamma_i^r\sim N(0,D_r^2)$ and $\gamma_i^s\sim N(0,D_s^2)$, respectively. In addition there is additive process noise $\xi_i^v\sim N(0,D_v^2)$, $\xi_i^q\sim N(0,D_q^2)$, and $\xi_i^{q^{\rm a}}\sim N(0,D_{q^{\rm a}}^2)$, $i=0,\ldots,N$, affecting the speed and flow equations, namely, (\ref{average speed}), and (\ref{flowtotal}), (\ref{flowvacs}), respectively. The noise statistics are summarized in Table \ref{table2}. 
 \begin{table}[t]
\caption{The measurement noise $\gamma_i^w$ and the process noise $\xi_i^w$, $i=0,\ldots,N$ affecting the $w$ variable at segment $i$. The variable $w$ can represent a flow (i.e., $w=q$, $w=r$, or $w=s$) or speed (i.e., $w=v$).}
\begin{center}
{\begin{tabular}{cc}
\hline\hline
Noise & Standard deviation\\
\hline
$\gamma_0^q$&$D_q=25\frac{\textrm{veh}}{\textrm{h}}$\\[2mm]
$\gamma_N^q$&$D_q=25\frac{\textrm{veh}}{\textrm{h}}$\\[2mm]
$\gamma_i^r$&$D_r=10\frac{\textrm{veh}}{\textrm{h}}$ \\[2mm]
$\gamma_i^s$ &$D_s=5\frac{\textrm{veh}}{\textrm{h}}$ \\[2mm]
 $\xi_i^v$&$D_v=5\frac{\textrm{km}}{\textrm{h}}$\\[2mm]
 $\xi_i^q$&$D_q=25\frac{\textrm{veh}}{\textrm{h}}$\\[2mm]
 $\xi_i^{q^{\rm a}}$&$D_{q^{\rm a}}=15\frac{\textrm{veh}}{\textrm{h}}$\\[2mm]
\hline
\end{tabular}}
\label{table2}
\end{center}
\end{table}%

\begin{figure}
\centering
\includegraphics[width=\linewidth]{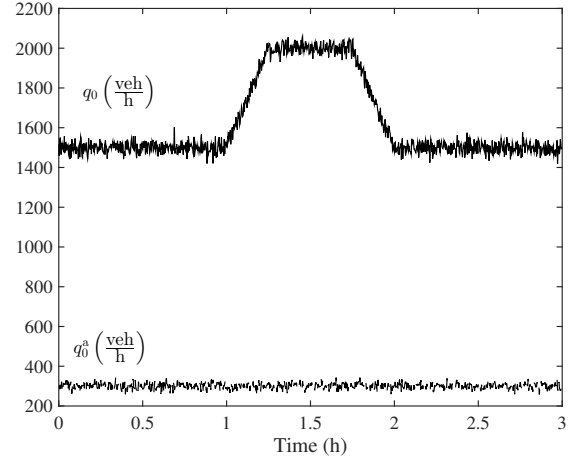}
\caption{The total flow of vehicles $q_0$ $\left(\mbox{in $\frac{\textrm{veh}}{\textrm{h}}$}\right)$ at the entry of the highway and the flow of connected vehicles $q^{\rm a}_0$ $\left(\mbox{in $\frac{\textrm{veh}}{\textrm{h}}$}\right)$ at the entry of the highway.}
\label{fig2}
\end{figure}

\begin{figure}
\centering
\includegraphics[width=\linewidth]{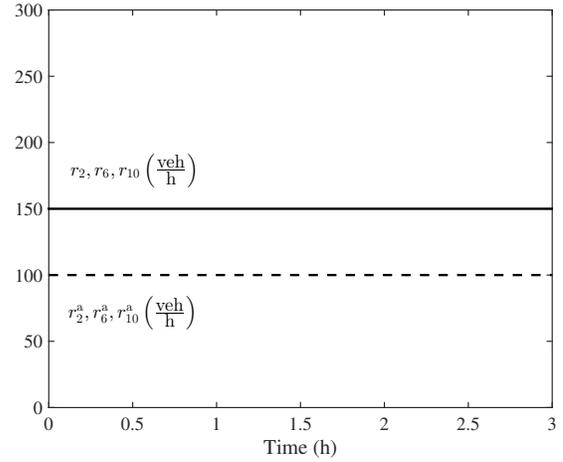}
\caption{The total flows $r_2, r_6, r_{10}$ $\left(\mbox{in $\frac{\textrm{veh}}{\textrm{h}}$}\right)$ of vehicles and the flows $r^{\rm a}_2,r^{\rm a}_6,r^{\rm a}_{10}$ $\left(\mbox{in $\frac{\textrm{veh}}{\textrm{h}}$}\right)$ of connected vehicles at the on-ramps.}
\label{fig3}
\end{figure}

The parameters and initial conditions of the Kalman filter (\ref{123})--(\ref{1234}), (\ref{adef})--(\ref{16}) are shown in Table \ref{table3}.
\begin{table}[t]
\caption{Parameters of the Kalman filter (\ref{123})--(\ref{1234}) and (\ref{adef})--(\ref{16}).}
\begin{center}
{\begin{tabular}{cc}
\hline\hline
Filter's parameter & Value\\
\hline
$Q$&$I_{N\times N}$\\[2mm]
$R$&$100$\\[2mm]
$\mu$&$(10,\ldots,10)^T$ \\[2mm]
$H$ &$I_{N\times N}$ \\[2mm]
\hline
\end{tabular}}
\label{table3}
\end{center}
\end{table}%
In Fig. \ref{fig2} we show the emptied scenario of input flow of connected vehicles and total input flow at the entry of the highway for our simulation investigation. We assume that there are three on-ramps at segments $2,6,10$. The total flow and the flow of connected vehicles at the on-ramps are shown in Fig. \ref{fig3}. Four off-ramps are supposedly present on the highway under study, specifically at segments $4,8,12$. It is assumed that $s_i=0.1 q_{i-1}$ and $s_i^{\rm a}=0.1q_{i-1}^{\rm a}$, $i=4,8,12$. 

\begin{figure}
\centering
\includegraphics[width=\linewidth]{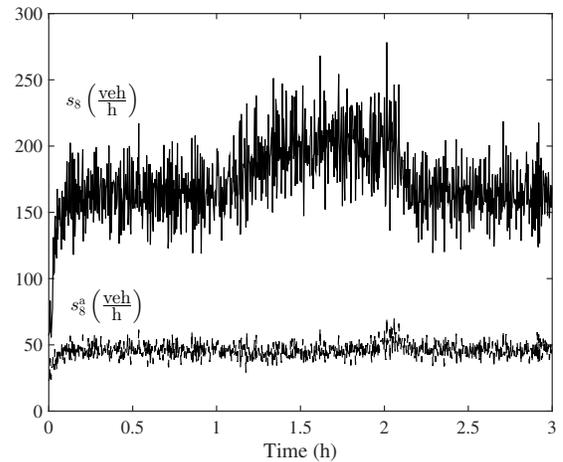}
\caption{The total flow $s_8$ $\left(\mbox{in $\frac{\textrm{veh}}{\textrm{h}}$}\right)$ of vehicles and the flow $s^{\rm a}_8$ $\left(\mbox{in $\frac{\textrm{veh}}{\textrm{h}}$}\right)$ of connected vehicles at the off-ramp located at the eighth segment.}
\label{fig4}
\end{figure}

\begin{figure}
\centering
\includegraphics[width=\linewidth]{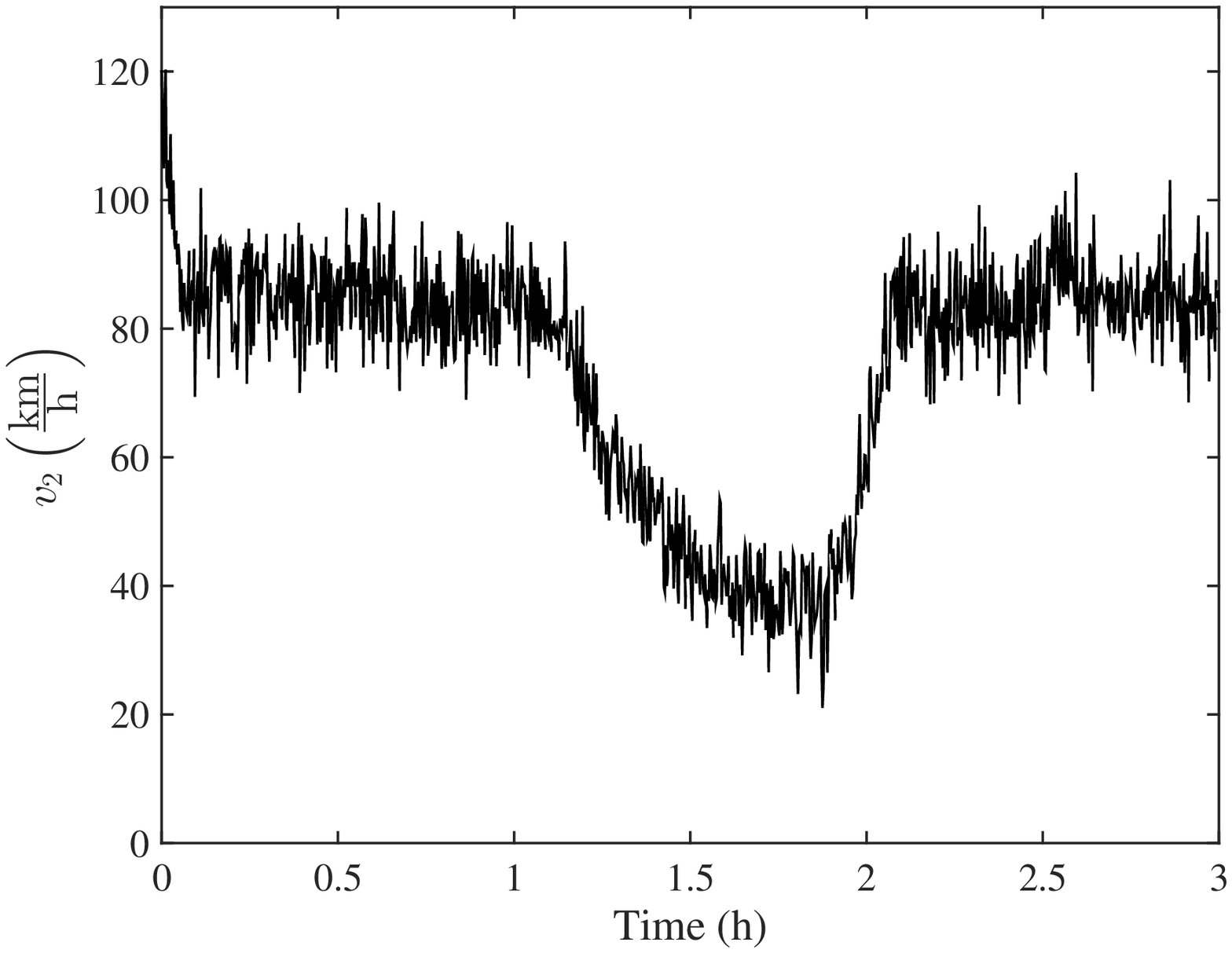}
\caption{The average speed $v_2$ $\left(\mbox{in $\frac{\textrm{km}}{\textrm{h}}$}\right)$ of the second segment of the highway as it is produced by the METANET model (\ref{eqrho}), (\ref{flowtotal}), (\ref{average speed}), (\ref{fund}) with parameters given in Table \ref{table1} and additive process noise given in Table \ref{table2}.}
\label{fig5}
\end{figure}

\begin{figure}
\centering
\includegraphics[width=\linewidth]{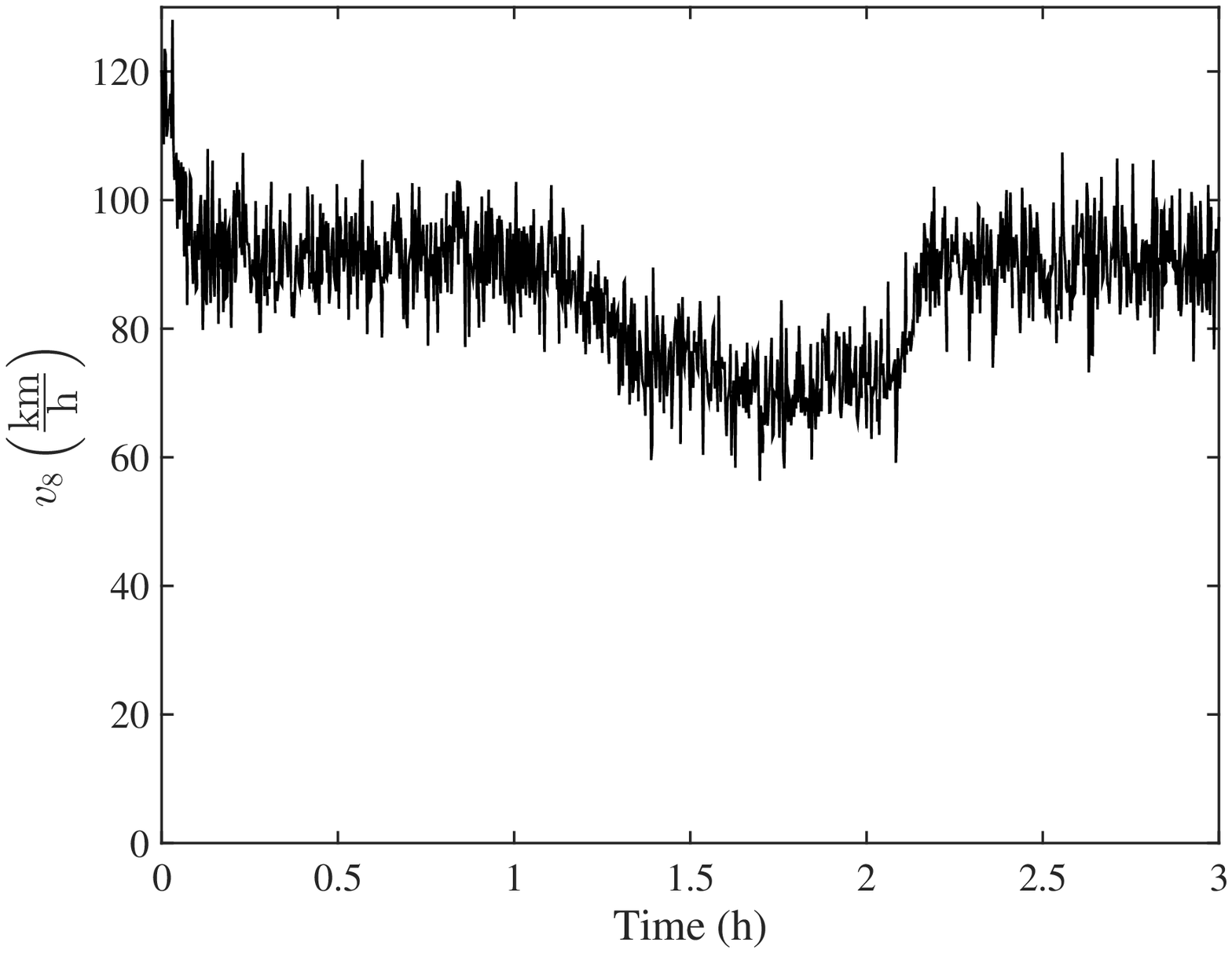}
\caption{The average speed $v_8$ $\left(\mbox{in $\frac{\textrm{km}}{\textrm{h}}$}\right)$ of the eighth segment of the highway as it is produced by the METANET model (\ref{eqrho}), (\ref{flowtotal}), (\ref{average speed}), (\ref{fund}) with parameters given in Table \ref{table1} and additive process noise given in Table \ref{table2}.}
\label{fig6}
\end{figure}

\begin{figure}
\centering
\includegraphics[width=\linewidth]{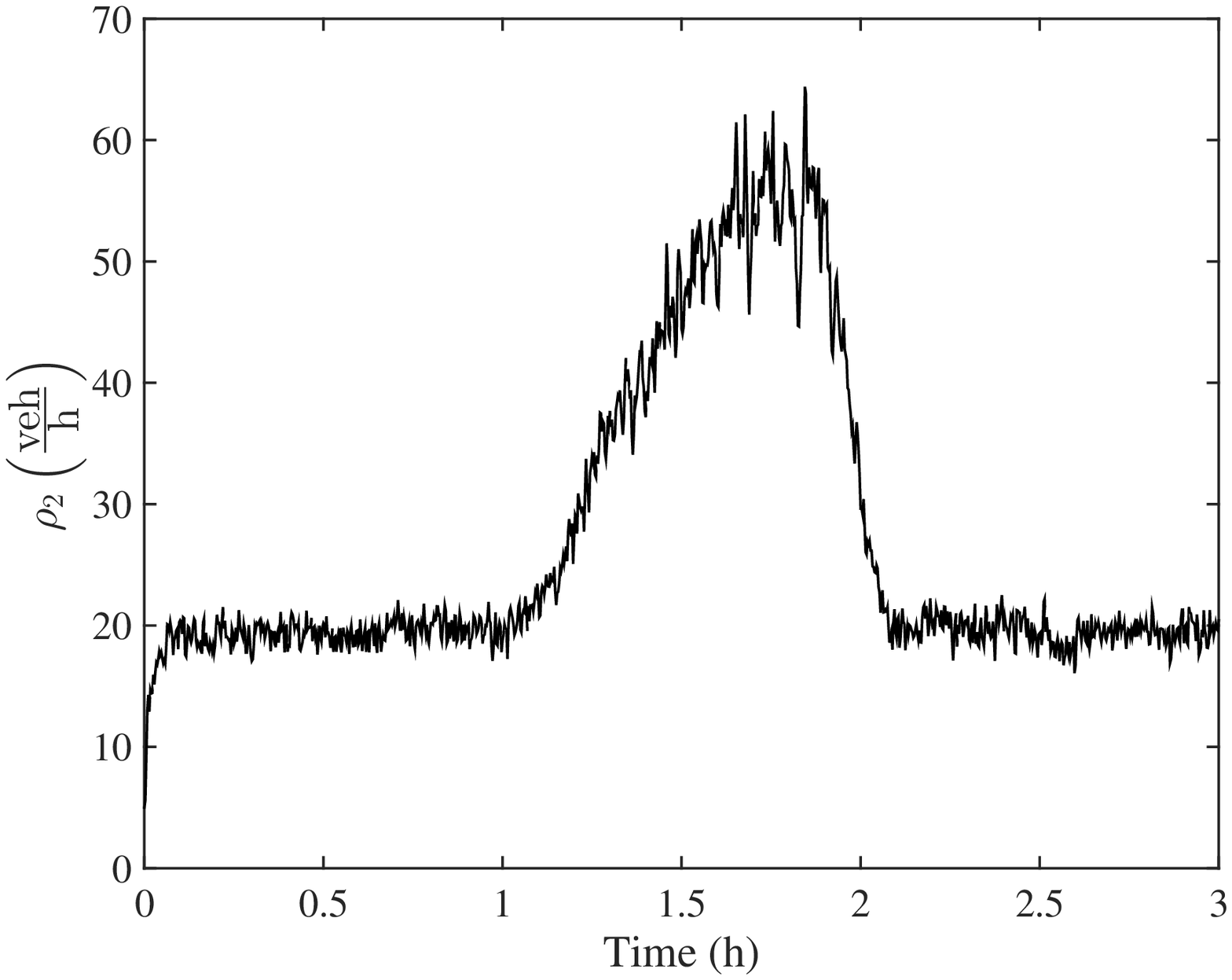}
\caption{The total density of vehicles $\rho_2$ $\left(\mbox{in $\frac{\textrm{veh}}{\textrm{km}}$}\right)$ at the second segment of the highway as it is produced by the METANET model (\ref{eqrho}), (\ref{flowtotal}), (\ref{average speed}), (\ref{fund}) with parameters given in Table \ref{table1} and additive process noise given in Table \ref{table2}.}
\label{fig7}
\end{figure}

\begin{figure}
\centering
\includegraphics[width=\linewidth]{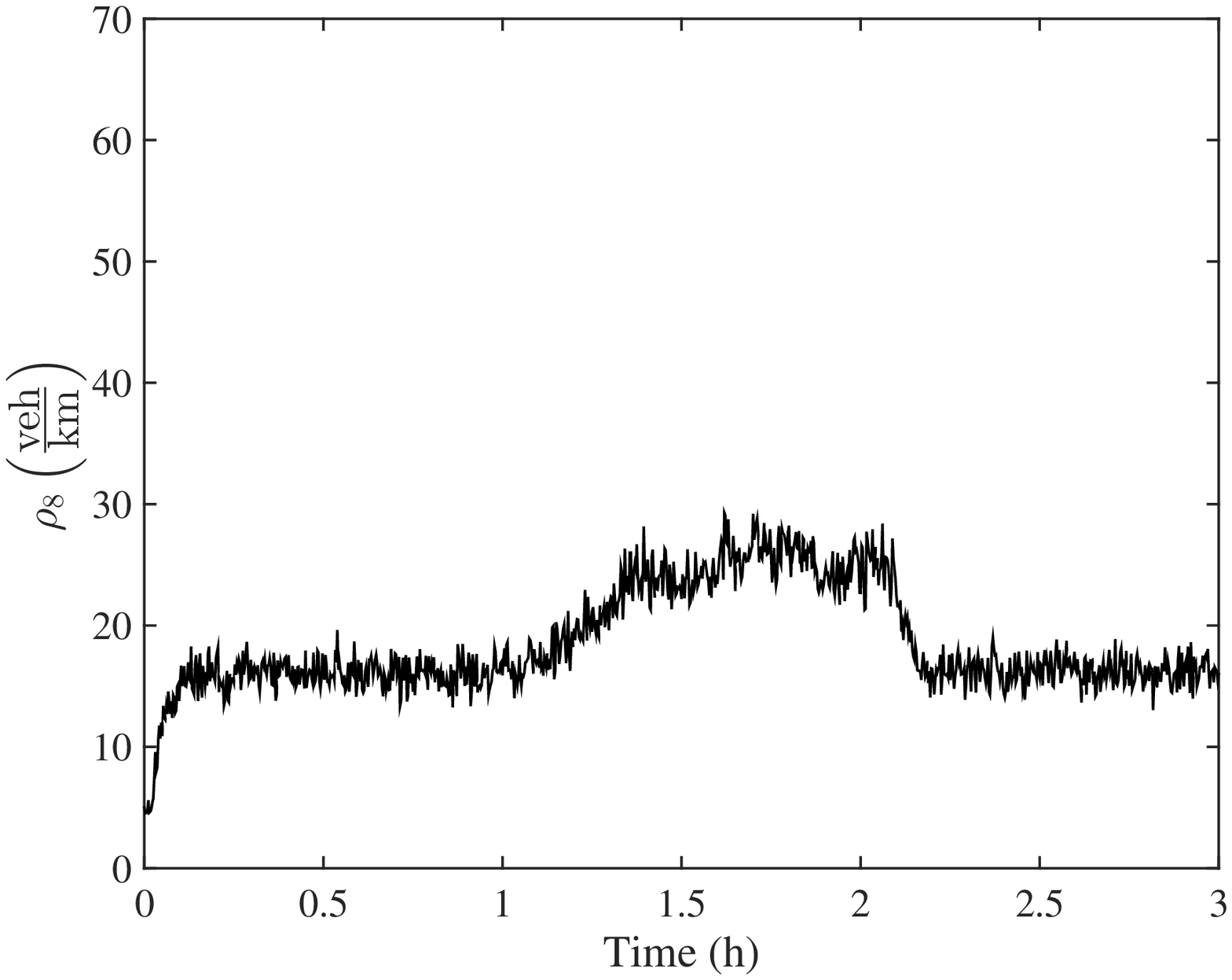}
\caption{The total density of vehicles $\rho_8$ $\left(\mbox{in $\frac{\textrm{veh}}{\textrm{km}}$}\right)$ at the eighth segment of the highway as it is produced by the METANET model (\ref{eqrho}), (\ref{flowtotal}), (\ref{average speed}), (\ref{fund}) with parameters given in Table \ref{table1} and additive process noise given in Table \ref{table2}.}
\label{fig8}
\end{figure}

The total flow and the flow of connected vehicles resulting from the simulation at the eighth off-ramp are shown in Fig. \ref{fig4}. The average speed at segments $2$ (where the first on-ramp is located) and $8$ (where the second off-ramp is located) are shown in Fig. \ref{fig5} and Fig. \ref{fig6}, respectively. The corresponding densities of the total number of vehicles are shown in Fig. \ref{fig7} and Fig. \ref{fig8}, respectively. It is evident from Fig. \ref{fig5} and Fig. \ref{fig7} that a congestion is created between the first and second hour of our test, whereas, free-flow conditions are reported for the first and last hour. Congestion starts approximately at the location of the second on-ramp, i.e., at the sixth segment of the highway, and propagates backwards all the way to the input of the highway.

\begin{figure}
\centering
\includegraphics[width=\linewidth]{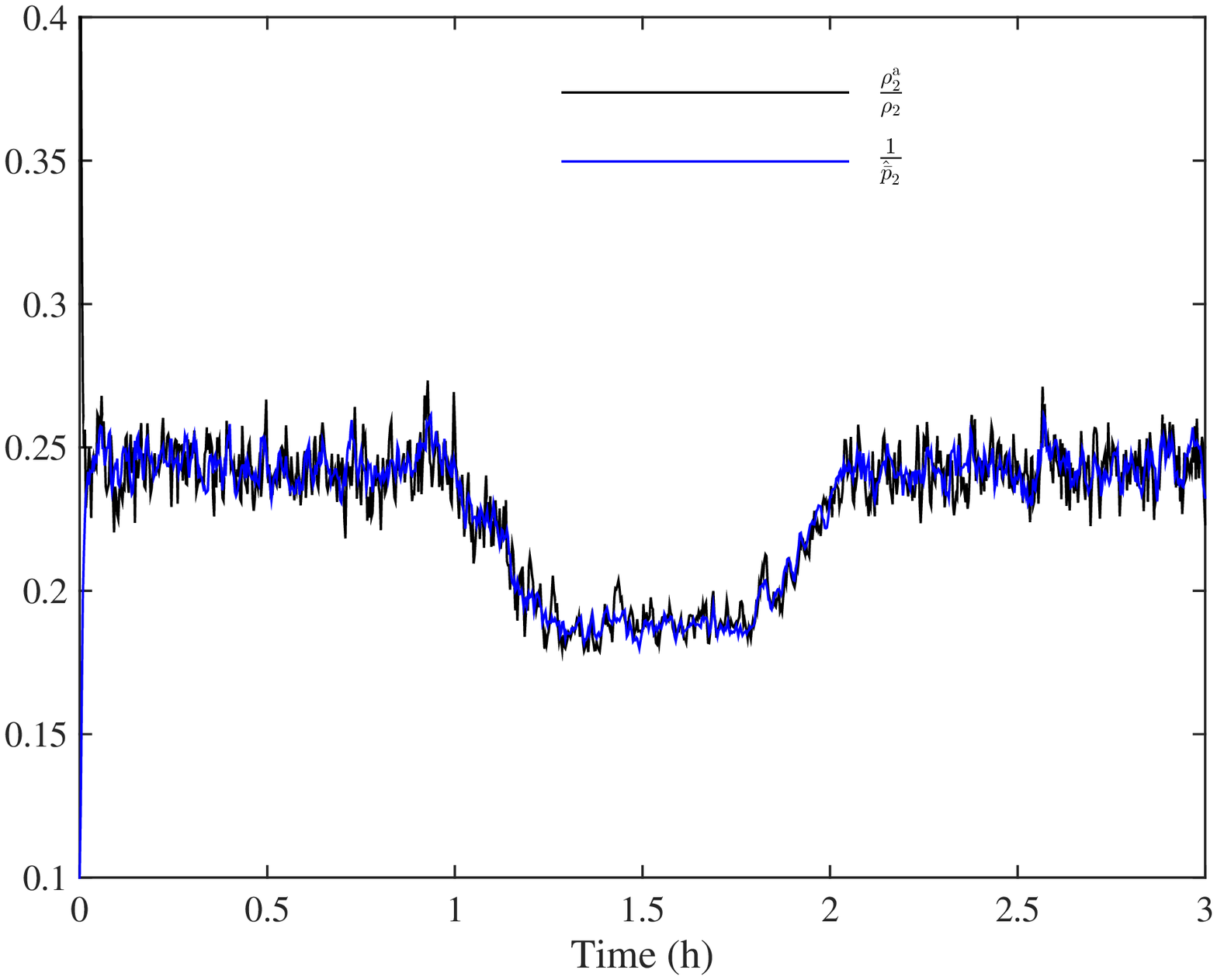}
\caption{The percentage of connected vehicles $\frac{\rho^{\rm a}_2}{\rho_2}$ on the second segment of the highway (black line) and its estimate $\frac{1}{\hat{\bar{p}}_2}$ (blue line) as it is produced by the Kalman filter with parameters given in Table \ref{table3}.}
\label{fig9}
\end{figure}

\begin{figure}
\centering
\includegraphics[width=\linewidth]{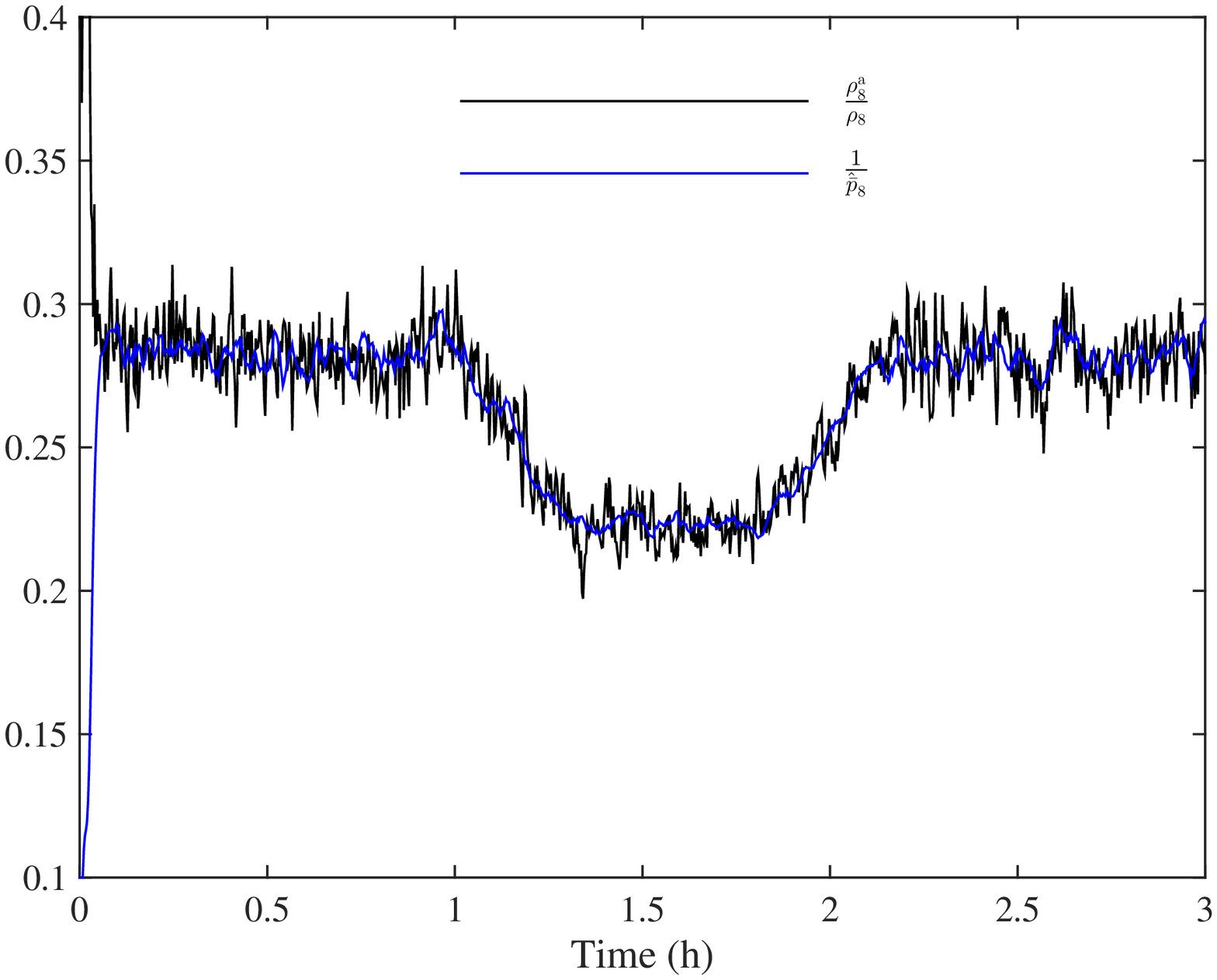}
\caption{The percentage of connected vehicles $\frac{\rho^{\rm a}_8}{\rho_8}$ on the eighth segment of the highway (black line) and its estimate $\frac{1}{\hat{\bar{p}}_8}$ (blue line) as it is produced by the Kalman filter with parameters given in Table \ref{table3}.}
\label{fig1nn}
\end{figure}

\begin{figure}
\centering
\includegraphics[width=\linewidth]{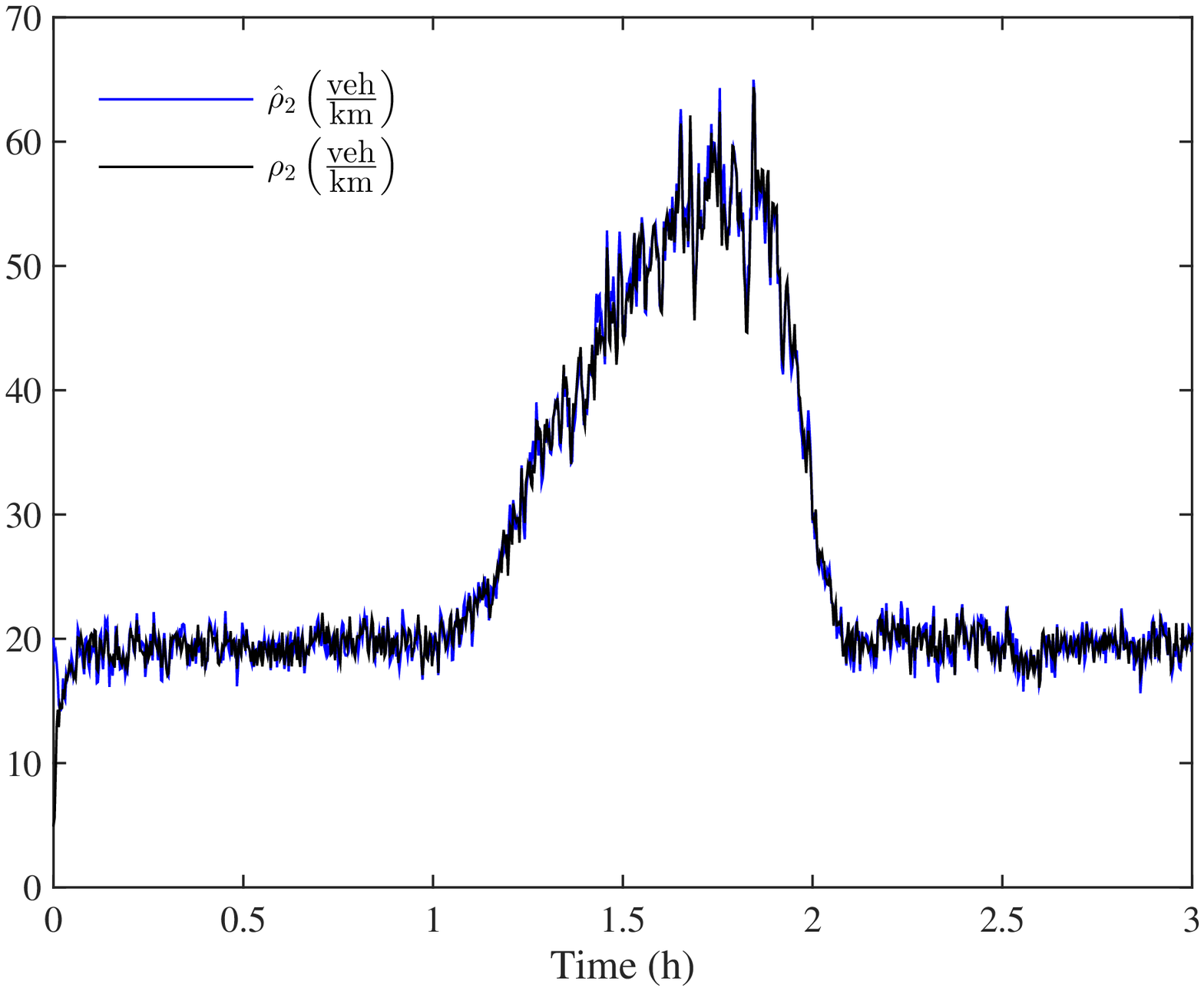}
\caption{The total density of vehicles $\rho_2$ $\left(\mbox{in $\frac{\textrm{veh}}{\textrm{km}}$}\right)$ on the second segment of the highway (black line) and its estimate $\hat{\rho}_2=\rho_2^{\rm a}\hat{\bar{p}}_2$ $\left(\mbox{in $\frac{\textrm{veh}}{\textrm{km}}$}\right)$ (blue line) as it is produced by the Kalman filter (\ref{123})--(\ref{1234}), (\ref{adef})--(\ref{16}) with parameters given in Table \ref{table3}.}
\label{fig1rhonn}
\end{figure}

\begin{figure}
\centering
\includegraphics[width=\linewidth]{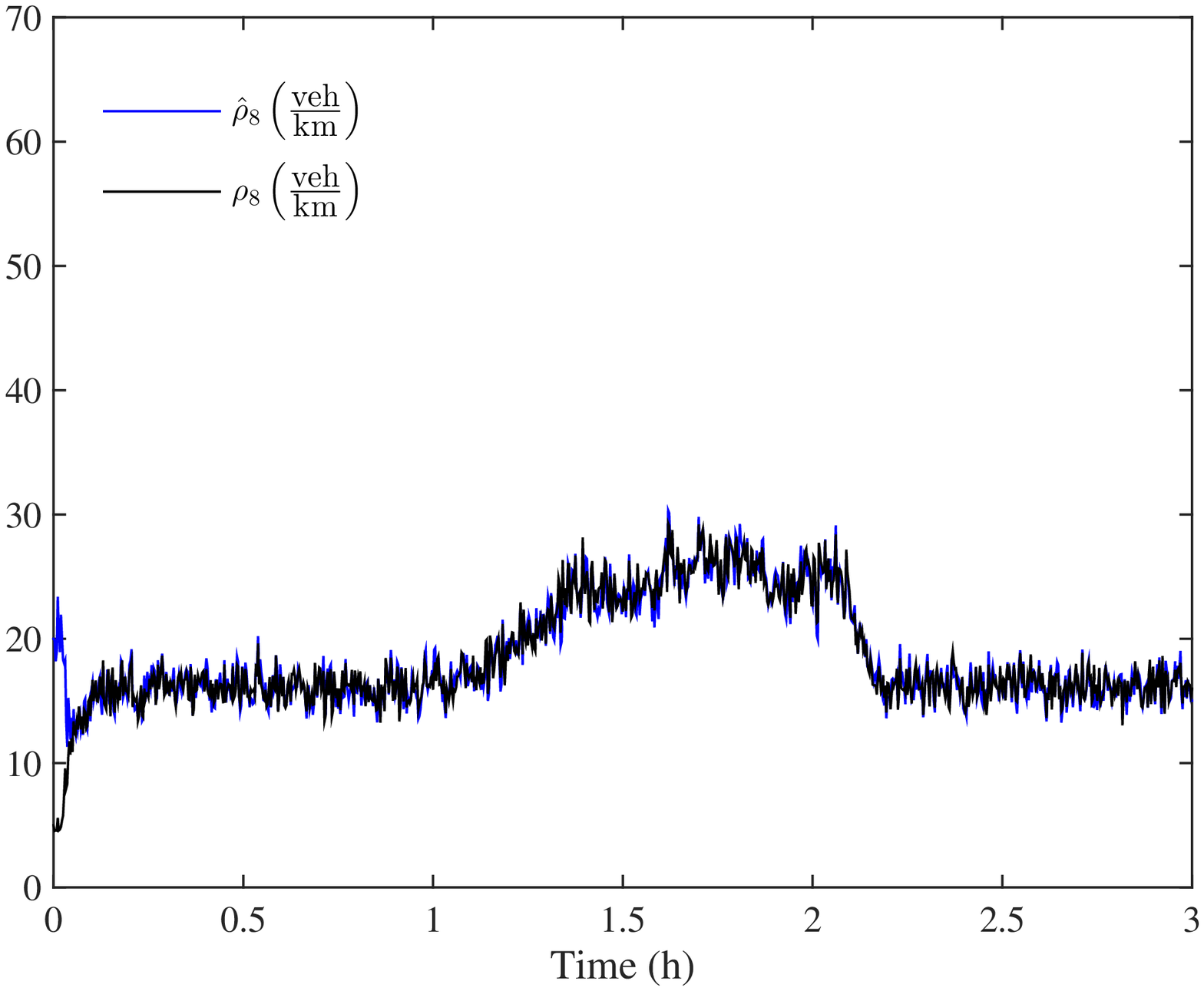}
\caption{The total density of vehicles $\rho_8$ $\left(\mbox{in $\frac{\textrm{veh}}{\textrm{km}}$}\right)$ on the eighth segment of the highway (black line) and its estimate $\hat{\rho}_8=\rho_8^{\rm a}\hat{\bar{p}}_8$ $\left(\mbox{in $\frac{\textrm{veh}}{\textrm{km}}$}\right)$ (blue line) as it is produced by the Kalman filter (\ref{123})--(\ref{1234}), (\ref{adef})--(\ref{16}) with parameters given in Table \ref{table3}.}
\label{fig1rho1nn}
\end{figure}

\begin{figure}
\centering
\includegraphics[width=\linewidth]{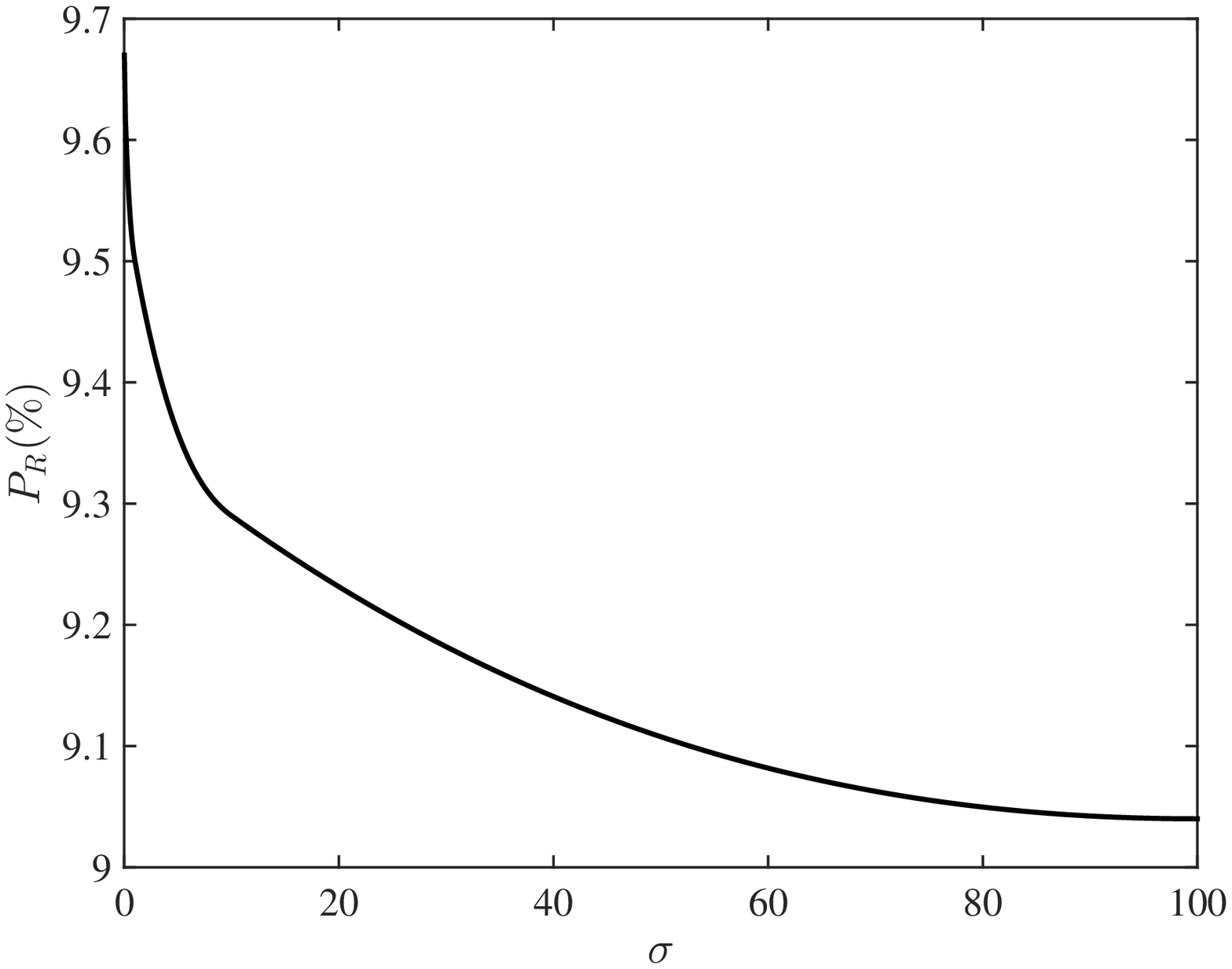}
\caption{The relative performance index $P_R$ (in \%) defined in (\ref{index}) as a function of the tuning parameter $Q=\sigma I_{N\times N}$ of the Kalman filter (\ref{123})--(\ref{1234}), (\ref{adef})--(\ref{16}) with parameters given in Table \ref{table3}.}
\label{perf}
\end{figure}

In both traffic conditions, our estimator successfully estimates the percentage of connected vehicles on the highway, as it is evident from Fig. \ref{fig9} and Fig. \ref{fig1nn}, which display the actual percentage and its estimate at two different segments of the highway, namely at segments $2$ (at which congested conditions prevail for one hour) and $8$, respectively. Note the very fast convergence of the produced percentage estimates, starting from remote initial values. Fig. \ref{fig1rhonn} and Fig. \ref{fig1rho1nn} display the resulting estimation of the total density of vehicles at segments $2$ and $8$, respectively, using relation (\ref{aa}). Moreover, in Fig. \ref{perf} we show the relative performance index of the estimation scheme defined as
\begin{eqnarray}
P_{R}=\frac{\sqrt{\frac{1}{MN}\sum_{k=0}^{k=M}\sum_{i=1}^{i=N}\left(\rho_i(k)-\rho_i^{\rm a}(k)\hat{\bar{p}}_i(k)\right)^2}}{\frac{1}{MN}\sum_{k=0}^{k=M}\sum_{i=1}^{i=N}\rho_i(k)},\label{index}
\end{eqnarray}
with simulation time horizon $M=\frac{3}{T}=1080$, as a function of the parameter $Q=\sigma I_{N\times N}$ of the Kalman filter while $R$ was kept constant at a value $R=100$. From Fig. \ref{perf} it is evident that the Kalman filter is robust to the choice of the tuning parameter $Q$.

\subsection{The Case of Unmeasured Total Flow at Off-Ramps}
\label{uns}
In the case that the total flow at off-ramps is not directly measured, one can extract this information as follows. We assume the following relations for the flow at off-ramps
\begin{eqnarray}
s_i&=&\beta_iq_{i-1}\label{s1}\\
s^{\rm a}_i&=&\beta^{\rm a}_iq^{\rm a}_{i-1},\label{s2}
\end{eqnarray}
where $\beta_i$ and ${\beta^{\rm a}_i}$ denote exit rates, i.e., the flow percentage of vehicles and connected vehicles, respectively, exiting at an off-ramp located at segment $i$. Assuming that the exit rates $\beta_i$ and $\beta^{\rm a}_i$ in (\ref{s1}), (\ref{s2}) are equal, and using (\ref{aa}), we get that 
\begin{eqnarray}
s_i=\beta^{\rm a}_iq^{\rm a}_{i-1}\bar{p}_{i-1}.\label{s3}
\end{eqnarray}
Substituting (\ref{s3}) into (\ref{per1}), we get a new model for $\bar{p}_i$, $i=1,\ldots,N$, of the form (\ref{12}), (\ref{120}), with ${C}$ given in (\ref{16}) and 
\setlength{\arraycolsep}{1pt}\begin{eqnarray}
{A}(k)&=&\left\{\begin{array}{lll}{a}_{ij}=\frac{T}{\Delta_i}\frac{\left(1-\beta^{\rm a}_i(k)\right)q^{\rm a}_{i-1}(k)}{\bar{g}_i^{\rm a}(k)},&\mbox{if $i-j=1$}\\&\mbox{and $i\geq2$}\\{a}_{ij}=\frac{\rho^{\rm a}_i(k)-\frac{T}{\Delta_i}q^{\rm a}_i(k)}{\bar{g}_i^{\rm a}(k)},&\mbox{if $i=j$}\\{a}_{ij}=0,&\mbox{otherwise}\end{array}\right\}\label{adef1}\\
{B}(k)&=&\left\{\begin{array}{lll}{b}_{ij}=\frac{T}{\Delta_i}\frac{1}{\bar{g}_1^{\rm a}(k)},&\mbox{if $i=1$ and $j=1,2$}\\{b}_{ij}=\frac{T}{\Delta_i}\frac{1}{\bar{g}_i^{\rm a}(k)},&\mbox{if $j-i=1$}\\{b}_{ij}=0,&\mbox{otherwise}\end{array}\right\}\label{defbu}\\
{u}(k)&=&\left[\begin{array}{cccc}q_0(k)&r_1(k)&\cdots&r_N(k)\end{array}\right]^T,\label{newu}
\end{eqnarray}\setlength{\arraycolsep}{5pt}where $\bar{g}_i^{\rm a}(k)=\rho_i^{\rm a}(k)+\frac{T}{\Delta_i}\left(\left(1-\beta_i^{\rm a}(k)\right)q_{i-1}^{\rm a}(k)-q_i^{\rm a}(k)\right)+\frac{T}{\Delta_i}r_i^{\rm a}(k)$, ${A}\in\mathbb{R}^{N\times N}$, and ${B}\in\mathbb{R}^{N\times (N+1)}$. With the same arguments as in Section \ref{seccal}, one can show that system (\ref{12}), (\ref{120}) with ${C}$ given in (\ref{16}) and ${A}$, ${B}$, and ${u}$ given in (\ref{adef1}), (\ref{defbu}), and (\ref{newu}), respectively, is observable provided that 
\begin{itemize}
\item relations (\ref{s1}) and (\ref{s2}) hold and 
\item $\beta_i^{\rm a}$ is available.
\end{itemize}
One can then implement the Kalman filter (\ref{123})--(\ref{1234}) with ${A}$, ${B}$, and ${u}$ given by (\ref{adef1}), (\ref{defbu}), and (\ref{newu}), respectively.

\section{Conclusions}
\label{sec conclusions}
We presented a macroscopic model-based approach for the estimation of the traffic state on highways in presence of connected vehicles, through estimating the percentage of connected vehicles, with respect to the total number of vehicles, on the highway. Specifically, we derived a linear time-varying system for the dynamics of the percentage and employed a Kalman filter for its estimation. We illustrated the effectiveness of our estimation design in simulation, using a second-order macroscopic traffic flow model as ground truth for the traffic state. We also discussed the possibility of applying our methodology to the case of unmeasured total flow of vehicles at off-ramps. 

A topic of ongoing research is the development of an alternative estimation algorithm for the total traffic density in highways utilizing only average speed measurements reported by connected vehicles, thus relaxing the requirement of measuring flows and densities for connected vehicles. That approach exploits the fact that the dynamics of the total density, as described by the conservation law equation, can be described by a linear time-varying system with known parameters on the basis of the same assumption employed in the present paper, namely, that the average speed of conventional vehicles is roughly equal to the average speed of connected vehicle. Our current research also includes several studies for the comparison of the performance of the two estimation approaches.

Future research will: 
\begin{itemize}
\item address the problem of unmeasured total flow of vehicles at on-ramps, via use of additional flow measurements at the mainstream of the highway,
\item address the problem of optimal fixed sensor placement on the highway;
\item validate the developed traffic estimation methodologies with a much more detailed microscopic simulation platform; considering a more realistic simulation of all involved real-time measurements.
\end{itemize}

\section*{Acknowledgments}
This research was supported by the European Research Council under the
European Union's Seventh Framework Programme (FP/2007-2013)/ERC Advanced Grant
Agreement n. 321132, project TRAMAN21.







\begin{thebibliography}{70}
\bibitem{alvarez}
L. Alvarez-lcaza, L. Munoz, X. Sun, and R. Horowitz, ``Adaptive observer for traffic density estimation," {\em ACC}, Boston, MA, 2004.
\bibitem{anger}
B. D. O. Anderson and J. B.  Moore, {\em Optimal Filtering}, Prentice-Hall, NJ, 1979.
\bibitem{atsarita}
V. Astarita, R. L. Bertini, S. d' Elia, and G. Guido, ``Motorway traffic parameter estimation from mobile phone counts," {\em European Journal of Operational Research}, vol. 175, pp. 1435--1446, 2006.
\bibitem{bose1}
A. Bose and P. Ioannou, ``Mixed manual/semi-automated traffic: a macroscopic analysis," {\em Transportation Research Part C}, vol. 11, pp. 439--462, 2003.
\bibitem{bose2}
A. Bose and P. Ioannou, ``Analysis of traffic flow with mixed manual and semiautomated vehicles," {\em IEEE Transactions on Intelligent Transportation Systems}, vol. 4, pp. 173--188, 2004.
\bibitem{davis1}
L. C. Davis, ``Effect of adaptive cruise control systems on mixed traffic flow near an on-ramp, {\em Physica A}, vol. 379, pp. 274--290, 2007.
\bibitem{fabri}
C. de Fabritiis, R. Ragona, and G. Valenti, ``Traffic estimation and prediction based on real time floating car data," {\em IEEE Conference on Intelligent Transportation Systems}, Beijing, China, 2008.
\bibitem{deng}
W. Deng, H. Lei, and X. Zhou, ``Traffic state estimation and uncertainty quantification based on heterogeneous data sources: A three detector approach," {\em Transp. Res. Part B}, vol. 57, pp. 132--157, 2013.
\bibitem{diakaki}
{C. Diakaki}, M. Papageorgiou, I. Papamichail, and I. K. Nikolos, ``Overview and analysis of vehicle automation and communication systems from a motorway traffic management perspective," {\em Transportation Research Part A}, to appear, 2015. Available at: http://www.traman21.tuc.gr/docs/wp1/TRAMAN21-D1-v4.pdf.
\bibitem{ge1}
J. I. Ge and G. Orosz, ``Dynamics of connected vehicle systems with delayed acceleration feedback," {\em Transportation Research Part C}, vol. 46, pp. 46--64, 2014.
\bibitem{bayen1}
J. C. Herrera and A. M. Bayen, ``Incorporation of Lagrangian measurements in freeway traffic state estimation," {\em Transportation Research Part B}, vol. 44, pp. 460--481, 2010.
\bibitem{bayen2}
J. C. Herrera, D. B. Work, R. Herring, X. Ban, Q. Jacobson, and A. M. Bayen, ``Evaluation of traffic data obtained via GPS-enabled mobile phones: The
Mobile Century field experiment," {\em Transportation Research Part C}, vol. 18, pp. 568--583, 2010.
\bibitem{heygi1}
A. Hegyi, D. Girimonte, R. Babuska, and B. De Schutter, ``A comparison of filter configurations for freeway traffic state estimation," {\em IEEE Conference on ITS}, Toronto, Canada, 2006.
\bibitem{kesting}
A. Kesting, M. Treiber, M. Schonhof, and D. Helbing, ``Adaptive cruise control design for active congestion avoidance," {\em Transportation Research Part C}, vol. 16, pp. 668--683, 2008.
\bibitem{miha}
L. Mihaylova, R. Boel, and A. Hegyi, ``Freeway traffic estimation within particle filtering framework," {\em Automatica}, vol. 43, pp. 290--300, 2007.
\bibitem{munoz}
L. Munoz, X. Sun, R. Horowitz, and L. Alvarez, ``Traffic density estimation with the Cell Transmission Model," {\em American Control Conference}, Denver, CO, 2003.
\bibitem{ngoduy4}
D. Ngoduy, S.P. Hoogendoorn, R. Liu, ``Continuum modeling of cooperative traffic flow dynamics," {\em Physica A: Statistical Mechanics and its Applications}, vol. 388, pp. 2705--2716, 2009.
\bibitem{ou}
Q. Ou, R. L. Bertini, J. W. C. van Lint, and S. P. Hoogendoorn, ``A theoretical framework for traffic speed estimation by fusing low-resolution probe vehicle data," {\em IEEE Transactions on Intelligent Transportation Systems}, vol. 12, pp. 747--756, 2011.
\bibitem{papageorge11}
M. Papageorgiou and A. Messmer, ``METANET: A macroscopic simulation program for motorway networks," {\em Traffic Engineering \& Control}, vol. 31, pp. 466--470, 1990.
\bibitem{rahmani1}
M. Rahmani, H. Koutsopoulos, and A. Ranganathan, ``Requirements and potential of GPS-based floating car data for traffic management: Stockholm case study," {\em IEEE Conference on Intelligent Transportation Systems}, Funchal, Portugal, 2010.
\bibitem{rajamani}
R. Rajamani and S. Shladover, ``An experimental comparative study of autonomous and cooperative vehicle-follower control systems," {\em Transportation Research Part C}, vol. 9, pp. 15--31, 2001.
\bibitem{rao}
B. Rao and P. Varaiya, ``Roadside intelligence for flow control in an intelligent vehicle and highway system," {\em Transportation Research Part C}, vol. 2, pp. 49--72, 1994.
\bibitem{roncoli0}
C. Roncoli, M. Papageorgiou, and I. Papamichail, ``Optimal control for multi-lane motorways in presence of Vehicle Automation and Communication Systems," {\em 19th IFAC World Congress}, Cape Town, South Africa, 2014.
\bibitem{roncoli00}
C. Roncoli, I. Papamichail, and M. Papageorgiou, ``Model predictive control for multi-lane motorways in presence of VACS," {\em IEEE Conference on Intelligent Transportation Systems}, Qingdao, China, 2014.
\bibitem{roncoli1}
C. Roncoli, M. Papageorgiou, and I. Papamichail, ``An optimisation-oriented first-order multi-lane model for motorway traffic," {\em Proceedings of the 94th Annual Meeting of the Transportation Research Board}, Washington, D.C.,  2015.
\bibitem{seo}
T. Seo, T. Kusakabe, and Y. Asakura, ``Estimation of flow and density using probe vehicles with spacing measurement equipment," {\em Transportation Research Part C}, vol. 53, pp. 134--150, 2015.
\bibitem{shladover}
S. E. Shladover, D. Su, and X.-Y. Lu, ``Impacts of cooperative adaptive cruise control on freeway traffic flow," {\em Transportation Research Record}, vol. 2324, pp. 63--70, 2012.
\bibitem{treiber}
M. Treiber, A. Kesting, and R. E. Wilson, ``Reconstructing the traffic state by fusion of heterogeneous data," {\em Computer-Aided Civil and Infrastructure Engineering}, vol. 26, pp. 408--419, 2011.
\bibitem{vanarem1}
B. van Arem, C. van Driel, and R. Visser, ``The impact of cooperative adaptive cruise control on traffic-flow characteristics," {\em IEEE Transactions on Intelligent Transportation Systems}, vol. 7, pp. 429--436, 2006.
\bibitem{varaya1}
P. Varaiya, ``Smart cars on smart roads: problems of control," {\em IEEE Transactions on Automatic Control}, vol. 38, 195--207, 1993.
\bibitem{wang van arem1}
M. Wang, W. Daamen, S. P. Hoogendoorn, and B. van Arem, ``Rolling horizon control framework for driver assistance systems. part II: cooperative sensing and cooperative control," {\em Transportation Research Part C}, vol. 40, pp. 290--311, 2014.
\bibitem{papageorge1}
Y. Wang and M. Papageorgiou, ``Real-time freeway traffic state estimation based on extended Kalman filter: a general approach," {\em Transportation Research Part B}, vol. 39, pp. 141--167, 2005.
\bibitem{papa new}
Y. Wang, M. Papageorgiou, A. Messmer, P. Coppola, A. Tzimitsi, and A. Nuzzolo, ``An adaptive freeway traffic state estimator," {\em Automatica}, vol. 45, pp. 10--24, 2009.
\bibitem{work1}
D. B. Work, O.-P. Tossavainen, S. Blandin, A. M. Bayen, T. Iwuchukwu, and K. Tracton, ``An ensemble Kalman filtering approach to highway traffic estimation using GPS enabled mobile devices," {\em IEEE Conference on Decision and Control}, Cancun, Mexico, 2008.
\bibitem{yi}
J. Yi and R. Horowitz, ``Macroscopic traffic flow propagation stability for adaptive cruise controlled vehicles," {\em Transportation Research Part C}, vol. 14, pp. 81--95, 2006.
\bibitem{yuan2}
Y. Yuan, J. W. C. van Lint, R. E. Wilson, F. van Wageningen-Kessels, and S. P. Hoogendoorn, ``Network-wide traffic state estimation using loop detector and floating car data," {\em Journal of Intelligent Transportation Systems}, vol. 18, pp. 41--50, 2014.
\end{thebibliography}
\end{document}